\documentclass[10pt, twoside]{amsart}
\usepackage{amsmath, amsthm, amscd, amsfonts, amssymb}
\usepackage{graphicx}

\begin{document}

\title{Fibonacci Statistical Convergence on Intuitionistic Fuzzy Normed Spaces}

\author{Murat Kiri\c{s}ci}

\address{[Murat Kiri\c{s}ci] Department of Mathematical Education, Hasan Ali Y\"{u}cel Education Faculty,
Istanbul University-Cerrahpa\c{s}a, Vefa, 34470, Fatih, Istanbul, Turkey \vskip 0.1cm }
\email{murat.kirisci@istanbul.edu.tr}

\begin{abstract}
In this paper, we study the concept of Fibonacci statistical convergence on intuitionisitic fuzzy normed space.
We define the Fibonacci statistically Cauchy sequences with respect to an intuitionisitic fuzzy normed space and
introduce the Fibonacci statistical completenes with respect to an intuitionisitic fuzzy normed space.
\end{abstract}

\maketitle

\pagestyle{plain} \makeatletter
\theoremstyle{plain}
\newtheorem{thm}{Theorem}[section]
\numberwithin{equation}{section}
\numberwithin{figure}{section}  %Comment out for sequentially-numbered
\theoremstyle{plain}
\newtheorem{pr}[thm]{Proposition}
\theoremstyle{plain}
\newtheorem{exmp}[thm]{Example}
\theoremstyle{plain}
\newtheorem{cor}[thm]{Corollary} %Delete [teo] to re-start numbering
\theoremstyle{plain}
\newtheorem{defin}[thm]{Definition}
\theoremstyle{plain}
\newtheorem{lem}[thm]{Lemma} %Delete [teo] to re-start numbering
\theoremstyle{plain}
\newtheorem{rem}[thm]{Remark}
\numberwithin{equation}{section}

% ----------------------------------------------------------------
\section{Introduction}

Fuzziness has revolutionized many areas such as mathematics, science, engineering, medicine. This concept was initiated by Zadeh \cite{Zadeh}. Not only Zadeh discovered  this concept, but he also developed the infrastructure of today's popular forms of use such as relations of similarity, decision making, and fuzzy programming in a short time.\\

Just like in the theory of sets, fuzzy sets(FS) also have led to the emergence of new mathematical concepts, research topics, and the design of engineering applications.
Therefore, the nature of the classical set theory must be well known and understood. In particular, consider the two fundamental laws of Boolean algebra the law of excluded middle and law of contradiction. In logic, the proposition every proposition is either true or false excludes any third, or middle, possibility, which gave this principle the name of  the law of excluded middle. When look at the principles of Boolean algebra, there are two items as prediction: "True" or "False". Whether classical, Boolean or crisp, set theory can be defined as a characteristic function of the membership of an  element $x$ in a set $A$.
For each elements of universal set $X$, the function that generates the values $0$ and $1$ is called the characteristic function.\\

In some real life problems in expert system, belief system, information fusion and so on , we must consider the truth-membership as well as the falsity-membership for proper description of an object in uncertain, ambiguous environment. Neither the fuzzy sets nor the interval valued fuzzy sets is appropriate for such a situation.
Intuitionistic fuzzy sets introduced by Atanassov \cite{Ata} is appropriate for such a situation. The intuitionistic fuzzy sets can only handle the incomplete
information considering both the truth-membership ( or simply membership ) and falsity-membership ( or non-membership ) values. It does not handle the indeterminate and inconsistent information which exists in belief system.\\

The notion of intuitionistic fuzzy metric space has been introduced by Park \cite{Park}. Furthermore, The concept of intuitionistic fuzzy normed space is given by Saadati and Park \cite{Saad}.

The study of statistical convergence was initiated by Fast\cite{Fast}. Schoenberg \cite{Sch} studied statistical convergence as
 a summability method and listed some of the elementary properties of statistical convergence. Both of these mathematicians
  mentioned that if a bounded sequence is statistically convergent to $L$ then it is Ces\`{a}ro summable to $L$. Statistical convergence
  also arises as an example of "convergence in density" as introduced by Buck\cite{Buck}. In \cite{Zygm}, Zygmund called this concept
  "almost convergence" and established relation between statistical convergence and strong summability. The idea of statistical convergence
  has been studied in different branches of mathematics such as number theory\cite{ErTe}, trigonometric series\cite{Zygm}, summability theory\cite{FreeSem},
  measure theory\cite{Mil}, Hausdorff locally convex topological vector spaces\cite{Maddox}. The concept of $\alpha\beta-$statistical convergence was introduced and studied by Aktu\v{g}lu\cite{Aktuglu}. In \cite{VatanAli}, Karakaya and Karaisa have been extended the concept of $\alpha\beta-$statistical convergence. Also,  they have been introduced the concept of weighted $\alpha\beta-$statistical convergence of order $\gamma$, weighted $\alpha\beta-$summmability of order of $\gamma$ and strongly weighted $\alpha\beta-$summable sequences of order $\gamma$, in \cite{VatanAli}. In \cite{Braha}, Braha gave a new weighted equi-statistical convergence and proved the Korovkin type theorems using the new definition.\\

In the present work, we will study Fibonacci statistical convergence, Fibonacci statistically Cauchy and Fibonacci statistical completeness
 with respect to the intuitionistic fuzzy normed space.\\

\section{Preliminaries}\label{chap:1}

The Fibonacci numbers are the sequence of numbers $(f_{n})$ for $n=1,2,\ldots$ defined by the linear recurrence equation
\begin{eqnarray*}
f_{n}=f_{n-1}+f_{n-2}  \quad n\geq 2,
\end{eqnarray*}

From this definition, it means that the first two numbers in Fibonacci sequence are either $1$ and $1$ (or $0$ and $1$)
 depending on the chosen starting point of the sequence and all subsequent number is the sum of the previous two. That is,
 we can choose $f_{1}=f_{2}=1$ or $f_{0}=0$, $f_{1}=1$. \\

Some of the fundamental properties of Fibonacci numbers are given as follows:
\begin{eqnarray*}
&&\lim_{n\rightarrow\infty}\frac{f_{n+1}}{f_{n}}=\frac{1+\sqrt{5}}{2}=\alpha, \quad \textrm{(golden ratio)}\\
&&\sum_{k=0}^{n}f_{k}=f_{n+2}-1  \quad  (n\in\mathbb{N}),\\
&&\sum_{k}\frac{1}{f_{k}} \quad \textrm{converges},\\
&&f_{n-1}f_{n+1}-f_{n}^{2}=(-1)^{n+1} \quad (n\geq 1) \quad \textrm{(Cassini formula)}
\end{eqnarray*}
It yields $f_{n-1}^{2}+f_{n}f_{n-1}-f_{n}^{2}=(-1)^{n+1}$, if we can substituting for $f_{n+1}$ in Cassini's formula.\\

Let $f_{n}$ be the $n$th Fibonacci number for every $n\in \mathbb{N}$. Then, we define the infinite matrix $\widehat{F}=(\widehat{f}_{nk})$ \cite{Kara1} by

\begin{eqnarray*}
\widehat{f}_{nk}=\left\{\begin{array}{ccl}
-\frac{f_{n+1}}{f_{n}}&, & (k=n-1)\\
\frac{f_{n}}{f_{n+1}}&, & (k=n)\\
0&, & (0\leq k < n-1 \textrm{or} k>n).
\end{array}\right.
\end{eqnarray*}

Let $A$ is a subset of positive integer. We consider the interval $[1,n]$ and select an integer in this interval, randomly. Then,
the ratio of the number of elements of $A$ in $[1,n]$ to the total number of elements in $[1,n]$ is belong to $A$, probably. For
$n\rightarrow \infty$, if this probability exists, that is this probability tends to some limit, then this limit is used to as
the asymptotic density of the set $A$. This mentions us that the asymptotic density is a kind of probability of choosing a number
from the set $A$.\\

Now, we give some definitions and properties of asymptotic density:\\

The set of positive integers will be denoted by $\mathbb{Z^{+}}$. Let $A$ and $B$ be subsets of $\mathbb{Z}^{+}$. If the symmetric
difference $A\Delta B$ is finite, then, we can say $A$ is asymptotically equal to $B$ and denote $A\sim B$. Freedman and Sember
have introduced the concept of a lower asymptotic density and defined a concept of convergence in density, in \cite{FreeSem}.
\begin{defin}\cite{FreeSem}
Let $f$ be a function which defined for all sets of natural numbers and take values in the interval $[0,1]$.
Then, the function $f$ is said to a lower asymptotic density, if the following conditions hold:
\begin{itemize}
\item[i.] $f(A)=f(B)$, if $A\sim B$,
\item[ii.] $f(A)+f(B)\leq f(A\cup B)$, if $A\cap B=\emptyset$,
\item[iii .] $f(A)+f(B)\leq 1+ f(A\cap B)$, for all $A$,
  \item[iv.]  $f(\mathbb{Z^{+}})=1$.
\end{itemize}
\end{defin}

We can define the upper density based on the definition of lower density as follows:\\

Let $f$ be any density. Then, for any set of natural numbers $A$, the function $\overline{f}$ is said to upper density associated with $f$, if
$\overline{f}(A)=1-f(\mathbb{Z}^{+} \backslash A)$.\\

Consider the set $A\subset \mathbb{Z}^{+}$. If $f(A)=\overline{f}(A)$, then we can say that the set $A$ has natural density
with respect to $f$. The term asymptotic density is often used for the function
\begin{eqnarray*}
d(A)=\liminf_{n\rightarrow\infty}\frac{A(n)}{n},
\end{eqnarray*}
where $A\subset \mathbb{N}$ and $A(n)=\sum_{a\leq n, a\in A}1$. Also the natural density of $A$ is given by $d(A)=\lim_{n}n^{-1}|A(n)|$ where
$|A(n)|$  denotes the number of elements in $A(n)$.\\

\begin{defin}
A real numbers sequence $x=(x_{k})$ is statistically convergent to $L$ provided that for every $\varepsilon >0$ the set
$\{n\in\mathbb{N}: |x_{n}-L|\geq \varepsilon\}$ has natural density zero. The set of all statistical convergent sequence is denoted by $S$.
In this case, we write $S-\lim x=L$ or $x_{k}\rightarrow L(S)$.\\
\end{defin}

\begin{defin}\cite{Fri}
The sequence $x=(x_{k})$ is statistically Cauchy sequence if for every  $\varepsilon >0$ there is a positive integer $N=N(\varepsilon)$ such
that
\begin{eqnarray*}
d\left(\{n\in\mathbb{N}: |x_{n}-x_{N(\varepsilon)}|\geq \varepsilon\}\right)=0.
\end{eqnarray*}

\end{defin}

It can be seen from the definition that statistical convergence is a generalization of the usual of notion of convergence that parallels
the usual theory of convergence. Kirisci and Karaisa \cite{KirKar} are defined Fibonacci statistical convergence and Fibonacci Cauchy sequence as follows:\\

\begin{defin}
A sequence $x=(x_{k})$ is said to be Fibonacci statistically
convergence(or $\widehat{F}-$statistically convergence) if there is
a  number $L$ such that for every $\epsilon> 0$ the set
$K_{\epsilon}(\widehat{F}):=\{k\leq
n:|\widehat{F}x_{k}-L|\geq\epsilon\}$ has natural density zero,
i.e., $d(K_{\epsilon}(\widehat{F}))=0$. That is

\begin{eqnarray*}
\lim_{n \to \infty}\frac{1}{n}\left|\{k\leq n:
|\widehat{F}x_{k}-L|\geq \epsilon \}\right|=0.
\end{eqnarray*}
\end{defin}

In this case we write $d(\widehat{F})-\lim x_{k}=L$ or $x_{k}\rightarrow L(S(\widehat{F}))$. The set of
$\widehat{F}-$statistically convergent sequences will be denoted by
$S(\widehat{F})$. In the case $L=0$, we will write
$S_{0}(\widehat{F})$.

\begin{defin}
Let $x=(x_{k})\in \omega$. The sequence $x$ is said to be $\widehat{F}-$statistically Cauchy if
there exists a number $N=N(\varepsilon)$ such that
\begin{eqnarray*}
\lim_{n \to \infty}\frac{1}{n}\left|\{k\leq n:
|\widehat{F}x_{k}-\widehat{F}x_{N}|\geq \epsilon \}\right|=0
\end{eqnarray*}
for every $\varepsilon >0$.\\
\end{defin}

In 1965, fuzzy set theory has long been presented to handle vagueness, inexact and imprecise data by Zadeh
\cite{Zadeh}. After Zadeh, fuzzy sets, especially fuzzy numbers have been widely studied and applied to various
fields, such as decision-making, pattern recognition, game theory and so on. In fuzzy sets, the degree
of memberships of the elements in a universe is a single value but, those single values cannot provide any
additional information because, in practice, information regarding elements corresponding to a fuzzy concept
may be incomplete. The fuzzy set theory is not capable of dealing with the lack of knowledge with respect
to degrees of memberships, Atanassov \cite{Ata} proposed the theory of intuitionistic fuzzy sets the extension of
Zadehs fuzzy sets by using a non-membership degree cope with the presence of vagueness and hesitancy
originating from imprecise knowledge or information.\\

\begin{defin}\cite{Schw}
Let $\bigstar$ be a binary operation with $\bigstar: [0,1]\times [0,1] \rightarrow [0,1]$ and given the following
properties:\\
\begin{itemize}
\item  $\forall x\in [0,1]$, $x \bigstar 1 = x$, \\
\item  $\forall x,y,z,t\in [0,1]$, $x\bigstar y\leq z\bigstar t$ whenever $x\leq y$ and $z\leq t$, \\
\item  $\bigstar$ is continuous,\\
\item  $\bigstar$ is commutative and associative.\\
\end{itemize}
If these conditions are satisfies, then  the operation $\bigstar$ is called a \emph{continuous t-norm}.\\
\end{defin}

\begin{defin}\cite{Schw}
Let $\blacktriangle$ be a binary operation with $\blacktriangle: [0,1]\times [0,1] \rightarrow [0,1]$ and given the following
properties:\\
\begin{itemize}
\item   $\forall x\in [0,1]$, $x \blacktriangle 0 = x$, \\
\item   $\forall x,y,z,t\in [0,1]$, $x\blacktriangle y\leq z\blacktriangle t$ whenever $x\leq y$ and $z\leq t$, \\
\item   $\blacktriangle$ is continuous,\\
\item   $\blacktriangle$ is commutative and associative.\\
\end{itemize}
If these conditions are satisfies, then  the operation $\blacktriangle$ is called a \emph{continuous t-conorm}.\\
\end{defin}

The concept of intuitionistic fuzzy normed space is given with the continuous t-norm and continuous t-conorm as follows:\\

\begin{defin}\cite{Saad}
Let $X$ be a vector space, $T$ and $U$ be two fuzzy sets on $X \times (0,\infty)$ and, $\bigstar$ and $\blacktriangle$
be a continuous t-norm and a continuous t-conorm, respectively. Consider the following conditions are hold:\\
\begin{itemize}
\item  $T(x,t)+U(x,t)\leq 1$,\\
\item  $T(x,t)>0$,\\
\item  $T(x,t)=1$ if and only if $x=0$,\\
\item  $T(\alpha x,t)=T\left(x, \frac{t}{|\alpha|}\right)$ for each $\alpha\neq 0$,\\
\item  $T(x,t)\bigstar T(y,s)\leq T(x+y,t+s)$,\\
\item  $T(x,.):(0,\infty)\rightarrow [0,1]$ is continuous,\\
\item  $\lim_{t\rightarrow\infty} T(x,t)=1$ and $\lim_{t\rightarrow 0} T(x,t)=0$,\\
\item  $U(x,t)<1$,\\
\item  $U(x,t)=0$ if and only if $x=0$,\\
\item  $U(\alpha x,t)=U\left(x, \frac{t}{|\alpha|}\right)$ for each $\alpha\neq 0$,\\
\item  $U(x,t)\blacktriangle U(y,s)\geq U(x+y,t+s)$,\\
\item  $U(x,.):(0,\infty)\rightarrow [0,1]$ is continuous,\\
\item  $\lim_{t\rightarrow\infty} U(x,t)=0$ and $\lim_{t\rightarrow 0} U(x,t)=1$.\\
\end{itemize}
for every $x,y\in X$ and $s,t>0$. Then, the 5-tuple $(X, T,S, \bigstar, \blacktriangle)$
is called \emph{intuitionistic fuzzy normed space}(IFNS). Therefore, $(T,U)$ is called an intuitionistic fuzzy norm(IFN).\\
\end{defin}

The concepts of convergence and Cauchy sequences in IFNS as follows  :\\
\begin{defin}\cite{Saad}
Consider the IFNS $(X, T,U, \bigstar, \blacktriangle)$. Let $\varepsilon, t>0$. If there exists $N\in \mathbb{N}$ such that $T(x_{k}-\mathcal{L}, t)>1-\varepsilon$
and $U(x_{k}-\mathcal{L},t)<\varepsilon$ for all $k\geq N$, then, in IFN $(T,U)$, a sequence $x=(x_{k})$ is said to be convergent to $\mathcal{L}\in X$.
It is denoted by $(T,U)-\lim x=\mathcal{L}$.\\

If there exists $N\in \mathbb{N}$ such that $T(x_{k}-x_{m}, t)>1-\varepsilon$
and $U(x_{k}-x_{m},t)<\varepsilon$ for all $k, m\geq N$,then, in IFN $(T,U)$, a sequence $x=(x_{k})$ is called Cauchy sequence.\\
\end{defin}

\section{Method}
In the theory of numbers, there are many different definitions of density. It is well known that the most popular
of these definitions is asymptotic density. But, asymptotic density does not exist for all sequences. New densities
have been defined to fill those gaps and to serve different purposes.\\

The asymptotic density is one of the possibilities to measure how large a subset of the set of natural number. We
know intuitively that positive integers are much more than perfect squares. Because, every perfect square is
positive and many other positive integers exist besides. However, the set of positive integers is not in fact larger
than the set of perfect squares: both sets are infinite and countable and can therefore be put in one-to-one correspondence.
Nevertheless if one goes through the natural numbers, the squares become increasingly scarce. It is precisely in this case,
natural density help us and makes this intuition precise.\\

The Fibonacci Sequence was firstly used in the Theory of Sequence Spaces by Kara and Ba\c{s}ar{\i}r\cite{KaraBas}.
Afterward, Kara\cite{Kara1} defined the Fibonacci
difference matrix $\widehat{F}$ by using the Fibonacci sequence $(f_{n})$ for $n\in \{0,1,\ldots\}$ and
introduced the new sequence spaces
related to the matrix domain of $\widehat{F}$. \\

The concept of the statistical convergence is studied  in the intuitionistic fuzzy normed space (\cite{KaSiMuFa}, \cite{KaSiMuFa2}, \cite{Karakus},
\cite{KoDe}, \cite{KuMur}, \cite{MoLo}, \cite{MuMo1}, \cite{MuMo2}). Kirisci and Karaisa \cite{KirKar} defined Fibonacci type statistical convergence and investigated
some fundamental properties. Further, in \cite{KirKar}, various approximation results concerning the classical
Korovkin theorem via Fibonacci type statistical convergence are given.\\

In this paper, by combining the definitions of Fibonacci sequence and statistical convergence with respect to the IFN $(T,U)$,
we will obtain a new concept of statistical convergence in IFNS and so on give a useful characterization for statistical convergence sequences on IFNS.
We will examine some basic properties of new statistical convergence in IFNS. We will see that this new statistical convergence method on IFNS is stronger than
the usual convergence on IFNS.\\

\section{Main Results}

\subsection{Fibonacci Statistical Convergence in IFNS}

Using the definitions and properties of density, statistical convergence, Fibonacci numbers and IFNS in Chapter 2 , we give new definitions as follows:\\

\begin{defin}\label{def:1}
Take an IFNS $(X,T,U,\bigstar,\blacktriangle)$. A sequence $x=(x_{k})$ is said to be Fibonacci statistical convergence with respect to IFN $(T,U)$ (briefly, FSC-IFN),
if there is a number $\mathcal{L}\in X$ such that for every $\varepsilon >0$ and $t>0$, the set
\begin{eqnarray*}
K_{\varepsilon}(\widehat{F}):=\big\{k\leq n: T(\widehat{F}x_{k}-\mathcal{L}, t)\leq 1-\varepsilon \quad \textit{or} \quad U(\widehat{F}x_{k}-\mathcal{L}, t)\geq\varepsilon\big\}
\end{eqnarray*}
has natural density zero, i.e., $d(K_{\varepsilon}(\widehat{F}))=0$. That is,
\begin{eqnarray*}
\lim_{n}\frac{1}{n} \bigg |\big\{k\leq n: T(\widehat{F}x_{k}-\mathcal{L}, t)\leq 1-\varepsilon  \quad \textit{or} \quad U(\widehat{F}x_{k}-\mathcal{L}, t)\geq\varepsilon\big\}  \bigg|=0.
\end{eqnarray*}
\end{defin}

In this case we write $d(\widehat{F})_{IFN}-\lim x_{k}=\mathcal{L}$ or $x_{k}\rightarrow \mathcal{L}(S(\widehat{F})_{IFN})$. The set of
FSC-IFN will be denoted by
$S(\widehat{F})_{IFN}$. In the case $\mathcal{L}=0$, we will write
$S_{0}(\widehat{F})_{IFN}$.\\

\begin{exmp}\label{exmp:1}
Let $(X,\|.\|)$ be a normed space and for all $a,b\in [0,1]$, $a\bigstar b=ab$ and $a\blacktriangle b=\min\{a+b,1\}$.
For all $x\in X$ and every $t>0$, consider
\begin{eqnarray*}
T(x,t):=\frac{t}{t+\|x\|}, \quad \quad U(x,t):=\frac{\|x\|}{t+\|x\|}.
\end{eqnarray*}
Then, $(X,T,U,\bigstar,\blacktriangle)$ is an IFNS. Define the $\widehat{F}x_{n}=(f_{n+1}^{2})=(1,2^{2},3^{2},5^{2},\cdots)$.
Since $f_{n+1}^{2}\rightarrow \infty$ as $k\rightarrow \infty$ and $\widehat{F}x=(1,0,0,\cdots)$, then $\widehat{F}x \in S$.
Consider
\begin{eqnarray*}
K_{n}(\varepsilon, t)=\{k\leq n: T_{F}(\widehat{F}x_{k},t)\leq 1-\varepsilon \quad \textit{or} \quad U_{F}(\widehat{F}x_{k},t)\geq \varepsilon\}
\end{eqnarray*}
for $\varepsilon \in (0,1)$ and for any $t>0$. When $n$ becomes sufficiently large, the quantity $T(\widehat{F}x_{k}-\mathcal{L}, t)$ becomes less than
$1-\varepsilon$ and similarly the quantity $U(\widehat{F}x_{k}-\mathcal{L}, t)$  becomes greater than $\varepsilon$. So, for $\varepsilon>0$ and $t>0$,
$K_{\varepsilon}(\widehat{F})=0$.
\end{exmp}

The following lemma can be easily get by using the definitions and properties of density given in Chapter 2 and Definition \ref{def:1}.\\

\begin{lem}\label{lem:1}
Take an IFNS $(X,T,U,\bigstar,\blacktriangle)$. For every $\varepsilon >0$ and $t>0$, the following statements are equivalent:
\begin{itemize}
\item[i.] $d(\widehat{F})_{IFN}-\lim x_{k}=\mathcal{L}$.
\item[ii.] $\lim_{n}\frac{1}{n} \bigg |\big\{k\leq n: T(\widehat{F}x_{k}-\mathcal{L}, t)\leq 1-\varepsilon\big\}\bigg|= \lim_{n}\frac{1}{n} \bigg |\big\{U(\widehat{F}x_{k}-\mathcal{L}, t)\geq\varepsilon\big\}  \bigg|=0$.
\item[iii.] $\lim_{n}\frac{1}{n} \bigg |\big\{k\leq n: T(\widehat{F}x_{k}-\mathcal{L}, t)> 1-\varepsilon\big\}\bigg|$  and $\lim_{n}\frac{1}{n} \bigg |\big\{k\leq n:U(\widehat{F}x_{k}-\mathcal{L}, t)<\varepsilon\big\}  \bigg|=1$.
\item[iv.] $\lim_{n}\frac{1}{n} \bigg |\big\{k\leq n: T(\widehat{F}x_{k}-\mathcal{L}, t)> 1-\varepsilon\big\}\bigg|= \lim_{n}\frac{1}{n} \bigg |\big\{U(\widehat{F}x_{k}-\mathcal{L}, t)<\varepsilon\big\}  \bigg|=1$.
\item[v.] $S-\lim T(\widehat{F}x_{k}-\mathcal{L}, t)=1$ and $S-\lim U(\widehat{F}x_{k}-\mathcal{L}, t)=0$.
\end{itemize}

\end{lem}

\begin{thm}
Take an IFNS $(X,T,U,\bigstar,\blacktriangle)$. The $d(\widehat{F})_{IFN}-\lim x_{k}=\mathcal{L}$ is unique,
when a sequence $x=(x_{k})$ is FSC-IFN.
\end{thm}

\begin{proof}
For $\mathcal{L}_{1}\neq \mathcal{L}_{2}$, suppose that $d(\widehat{F})_{IFN}-\lim x_{k}=\mathcal{L}_{1}$ and $d(\widehat{F})_{IFN}-\lim x_{k}=\mathcal{L}_{2}$. Choose $r>0$ such that $(1-r)\bigstar(1-r)>1-\varepsilon$ and $r\blacktriangle r<\varepsilon$, for a given $\varepsilon>0$.
Define the following set, for any $t>0$:\\
\begin{eqnarray*}
&&K_{T,1}(r,t):=\big\{k\in \mathbb{N}: T(\widehat{F}x_{k}-\mathcal{L}_{1}, \frac{t}{2})\leq 1-r\big\},\\
&&K_{T,2}(r,t):=\big\{k\in \mathbb{N}: T(\widehat{F}x_{k}-\mathcal{L}_{2}, \frac{t}{2})\leq 1-r\big\},\\
&&K_{U,1}(r,t):=\big\{k\in \mathbb{N}: U(\widehat{F}x_{k}-\mathcal{L}_{1}, \frac{t}{2})\geq r\big\},\\
&&K_{U,1}(r,t):=\big\{k\in \mathbb{N}: U(\widehat{F}x_{k}-\mathcal{L}_{2}, \frac{t}{2})\geq r\big\}.
\end{eqnarray*}
Then, for all $t>0$ and using the Lemma \ref{lem:1}, we have
\begin{eqnarray*}
d\big(K_{T,1}(\varepsilon,t))=d(K_{U,1}(\varepsilon,t)\big)=0
\end{eqnarray*}
because of $d(\widehat{F})_{IFN}-\lim x_{k}=\mathcal{L}_{1}$.
Similarly,  for all $t>0$, we have
\begin{eqnarray*}
d\big(K_{T,2}(\varepsilon,t))=d(K_{U,2}(\varepsilon,t)\big)=0
\end{eqnarray*}
because of $d(\widehat{F})_{IFN}-\lim x_{k}=\mathcal{L}_{2}$.
Let
\begin{eqnarray*}
K_{T,U}(\varepsilon,t):=\big\{K_{T,1}(\varepsilon,t)\cup K_{T,2}(\varepsilon,t)\big\}\cap\big\{K_{U,1}(\varepsilon,t)\cup K_{U,2}(\varepsilon,t)\big\}.
\end{eqnarray*}
Hence reveal that $d(K_{T,U}(\varepsilon,t))=0$ which implies $d(\mathbb{N} / K_{T,U}(\varepsilon,t))=1$.
Then, we have two possible cases, when take $k\in \mathbb{N} / K_{T,U}(\varepsilon,t)$:
\begin{itemize}
\item[i.] $k\in \mathbb{N} / \big(K_{T,1}(\varepsilon,t)\cup K_{T,2}(\varepsilon,t)\big)$,\\
\item[ii.] $k\in \mathbb{N} / \big(K_{U,1}(\varepsilon,t)\cup K_{U,2}(\varepsilon,t)\big)$.
\end{itemize}

Firstly, consider (i). Then, we have
\begin{eqnarray*}
T(\mathcal{L}_{1}-\mathcal{L}_{2},t)\geq T\left(\widehat{F}x_{k}-\mathcal{L}_{1},\frac{t}{2}\right)\bigstar
T\left(\widehat{F}x_{k}-\mathcal{L}_{2},\frac{t}{2}\right)> (1-r)\bigstar (1-r).
\end{eqnarray*}

In this case,
\begin{eqnarray}\label{equ:1}
T(\mathcal{L}_{1}-\mathcal{L}_{2}, t)> 1- \varepsilon,
\end{eqnarray}
because of $(1-r)\blacktriangle (1-r)>1- \varepsilon$.

Using the (\ref{equ:1}), for all $t>0$, we obtain $T(\mathcal{L}_{1}-\mathcal{L}_{2}, t)=1$, where $\varepsilon >0$ is arbitrary.
From here, $\mathcal{L}_{1}=\mathcal{L}_{2}$ is obtained.\\

Moreover, then we can write that
\begin{eqnarray*}
U(\mathcal{L}_{1}-\mathcal{L}_{2},t)\leq U\left(\widehat{F}x_{k}-\mathcal{L}_{1},\frac{t}{2}\right)\blacktriangle
U\left(\widehat{F}x_{k}-\mathcal{L}_{2},\frac{t}{2}\right)< r \blacktriangle r.
\end{eqnarray*}
if we take $k\in \mathbb{N} / \big(K_{U,1}(\varepsilon,t)\cup K_{U,2}(\varepsilon,t)\big)$. Using $r \blacktriangle r< \varepsilon$,
we see that
\begin{eqnarray*}
U(\mathcal{L}_{1}-\mathcal{L}_{2},t)<\varepsilon.
\end{eqnarray*}
For all $t>0$, we obtain $U(\mathcal{L}_{1}-\mathcal{L}_{2}, t)=0$, where $\varepsilon >0$ is arbitrary.
Thus, $\mathcal{L}_{1}=\mathcal{L}_{2}$ and this step completes the proof.\\

\end{proof}

\begin{thm}\label{theo:1}
Take an IFNS $(X,T,U,\bigstar,\blacktriangle)$. If  $(T,U)-\lim x=\mathcal{L}$, then
$d(\widehat{F})_{IFN}-\lim x_{k}=\mathcal{L}$.
\end{thm}

\begin{proof}
Let $(T,U)-\lim x=\mathcal{L}$. Then, for every $\varepsilon>0$ and $t>0$, there is a number $N\in \mathbb{N}$ such that
$T(x_{k}-\mathcal{L}, t)>1-\varepsilon$ and $U(x_{k}-\mathcal{L}, t)<\varepsilon$, for all $k\geq N$.
It shows that the set
\begin{eqnarray*}
\big\{k\in \mathbb{N}: T(x_{k}-\mathcal{L}, t)\leq 1-\varepsilon \quad \textit{or} \quad U(x_{k}-\mathcal{L}, t)\geq\varepsilon \big\}
\end{eqnarray*}
has at most finitely many terms. Then,
\begin{eqnarray*}
\lim_{n}\frac{1}{n} \bigg |\{k\leq n: T(\widehat{F}x_{k}-\mathcal{L}, t)\leq 1-\varepsilon  \quad \textit{or} \quad U(\widehat{F}x_{k}-\mathcal{L}, t)\geq\varepsilon\}  \bigg|=0,
\end{eqnarray*}
because of every finite subset of the natural numbers has density zero. This completes the proof.
\end{proof}

The converse of this theorem need not be true. Example \ref{exmp:1} shows that the converse of
Theorem \ref{theo:1} is not valid.

\begin{thm}\label{theo:2}
Let $(X,T,U,\bigstar,\blacktriangle)$ be an IFNS. $d(\widehat{F})_{IFN}-\lim x_{k}=\mathcal{L}$ if and only if
there exists a increasing index sequence $J=\{j_{1},j_{2},\cdots\}\subseteq \mathbb{N}$, when $d(J)=1$ and $(T,U)-\lim_{n\rightarrow\infty} x_{j_{n}}=\mathcal{L}$.
\end{thm}

\begin{proof}
Assume that $d(\widehat{F})_{IFN}-\lim x_{k}=\mathcal{L}$. For any $t>0$ and $u=1,2,\ldots$

\begin{eqnarray*}
M_{T,U}(u,t)=\bigg \{k\in \mathbb{N}: T(\widehat{F}x_{k}-\mathcal{L}, t)> 1-\frac{1}{u} \quad \textit{or} \quad U(\widehat{F}x_{k}-\mathcal{L}, t)< \frac{1}{u} \bigg\}.
\end{eqnarray*}

and

\begin{eqnarray*}
P_{T,U}(u,t)=\bigg \{k\in \mathbb{N}: T(\widehat{F}x_{k}-\mathcal{L}, t)\leq 1-\frac{1}{u} \quad \textit{or} \quad U(\widehat{F}x_{k}-\mathcal{L}, t)\geq \frac{1}{u} \bigg\}.
\end{eqnarray*}

Then, $d(P_{T,U}(u,t))=0$, since $d(\widehat{F})_{IFN}-\lim x_{k}=\mathcal{L}$. Further, for $t>0$ and $u=1,2,\ldots$,
\begin{eqnarray*}
M_{T,U}(u,t)\supset M_{T,U}(u+1,t).
\end{eqnarray*}

Since $d(\widehat{F})_{IFN}-\lim x_{k}=\mathcal{L}$, it is clear that, for $t>0$ and $u=1,2,\ldots$,
\begin{eqnarray}\label{equ:2}
d(M_{T,U}(u,t))=1.
\end{eqnarray}
Now, we will show that $k\in M_{T,U}(u,t)$, $x_{k}\rightarrow \mathcal{L}(S(\widehat{F})_{IFN})$.
Suppose that for some $k\in M_{T,U}(u,t)$, $x_{k}\nrightarrow \mathcal{L}(S(\widehat{F})_{IFN})$.
Therefore there is $\mu >0$ and a positive integer $N$ such that for all $k\geq N$,
$T(\widehat{F}x_{k}-\mathcal{L}, t)\leq 1-\mu$ or $U(\widehat{F}x_{k}-\mathcal{L}, t)\geq \mu$.
Let $T(\widehat{F}x_{k}-\mathcal{L}, t)> 1-\mu$ or $U(\widehat{F}x_{k}-\mathcal{L}, t)< \mu$, for all $k\geq N$.
Then, $d(\{k\in \mathbb{N}: T(\widehat{F}x_{k}-\mathcal{L}, t)> 1-\mu \quad \textit{and} \quad U(\widehat{F}x_{k}-\mathcal{L}, t)< \mu\})=0$.
Since $\mu>\frac{1}{u}$, we have $d(M_{T,U}(u,t))=0$, which contradicts (\ref{equ:2}). Therefore
$k\in M_{T,U}(u,t)$, $x_{k}\rightarrow \mathcal{L}(S(\widehat{F})_{IFN})$.\\

Suppose that there exists a subset $J=\{j_{1},j_{2},\cdots\}\subseteq \mathbb{N}$ such that $d(J)=1$
and $(T,U)-\lim_{n\rightarrow \infty}x_{j_{n}}=\mathcal{L}$, i.e. there exists $N\in \mathbb{N}$ such that
for every $\mu>0$ and $t>0$, $T(\widehat{F}x_{k}-\mathcal{L}, t)> 1-\mu$ or $U(\widehat{F}x_{k}-\mathcal{L}, t)< \mu$.
Now,
\begin{eqnarray*}
K_{T,U}(\mu, t):=\big \{k\in \mathbb{N}:  T(\widehat{F}x_{k}-\mathcal{L}, t)\leq 1-\mu \quad \textit{or} \quad U(\widehat{F}x_{k}-\mathcal{L}, t)\geq \mu \big\}\\
\subseteq \mathbb{N} - \{j_{N+1},j_{N+2},\ldots\}.
\end{eqnarray*}
Therefore $d(M_{T,U}(\mu, t))\leq 1-1=0$. Hence $d(\widehat{F})_{IFN}-\lim x_{k}=\mathcal{L}$.
\end{proof}

\subsection{Fibonacci Statistical Complete IFNS}

\begin{defin}\label{def:2}
Take an IFNS $(X,T,U,\bigstar,\blacktriangle)$. A sequence $x=(x_{k})$ is said to be Fibonacci statistical Cauchy with respect to IFN $(T,U)$ (briefly, $FSCa-IFN$),
if for every $\varepsilon>0$ and $t>0$, there exists $N=N(\varepsilon)$ such that
\begin{eqnarray*}
K_{\varepsilon}(\widehat{F}):=\{k\leq n: T(\widehat{F}x_{k}-\widehat{F}x_{N}, t)\leq 1-\varepsilon \quad \textit{or} \quad U(\widehat{F}x_{k}-\widehat{F}x_{N}, t)\geq\varepsilon\}
\end{eqnarray*}
has natural density zero, i.e., $d(K_{\varepsilon}(\widehat{F}))=0$. That is,
\begin{eqnarray*}
\lim_{n}\frac{1}{n} \bigg |\{k\leq n: T(\widehat{F}x_{k}-\widehat{F}x_{N}, t)\leq 1-\varepsilon  \quad \textit{or} \quad U(\widehat{F}x_{k}-\widehat{F}x_{N}, t)\geq\varepsilon\}  \bigg|=0.
\end{eqnarray*}
\end{defin}

\begin{thm}\label{theo:3}
Let $(X,T,U,\bigstar,\blacktriangle)$ be an IFNS. If a sequence $x=(x_{k})$ is $FSC-IFN$, then it is
$FSCa-IFN$.
\end{thm}

\begin{proof}
 Let $d(\widehat{F})_{IFN}-\lim x_{k}=\mathcal{L}$. For a given $\varepsilon>0$, choose $u>0$
 such that $(1-\varepsilon)\bigstar(1-\varepsilon)>1-u$ and $\varepsilon\blacktriangle\varepsilon<u$.
 Then, for $t>0$, we have,
 \begin{eqnarray}\label{equ:3}
d\left(A(\varepsilon, t)\right)=d(\{k\in \mathbb{N}: T(\widehat{F}x_{k}-\mathcal{L}, t/2))\leq 1-\varepsilon
 \quad \textit{or} \quad U(\widehat{F}x_{k}-\mathcal{L}, t/2)\geq \varepsilon\})=0
 \end{eqnarray}
 which implies that
 \begin{eqnarray*}
 d\left(A^{C}(\varepsilon, t)\right)=d(\{k\in \mathbb{N}: T(\widehat{F}x_{k}-\mathcal{L}, t/2))> 1-\varepsilon
 \quad \textit{and} \quad U(\widehat{F}x_{k}-\mathcal{L}, t/2)< \varepsilon\})=1.
 \end{eqnarray*}
 Let $q\in A^{C}(\varepsilon, t)$. Then
 \begin{eqnarray*}
 T(\widehat{F}x_{q}-\mathcal{L}, t)> 1-\varepsilon \quad \textit{and} \quad U(\widehat{F}x_{q}-\mathcal{L}, t)< \varepsilon.
 \end{eqnarray*}
 Now, let
 \begin{eqnarray*}
 B(\varepsilon, t)=\left\{ k\in \mathbb{N}:T(\widehat{F}x_{k}-\widehat{F}x_{q}, t)\leq 1-u
 \quad \textit{or} \quad U(\widehat{F}x_{k}-\widehat{F}x_{q}, t)\geq u\right\}.
 \end{eqnarray*}
We need to show that $B(\varepsilon, t)\subset A(\varepsilon, t)$. Let $k\in B(\varepsilon, t)/A(\varepsilon, t)$.
Then we have
\begin{eqnarray*}
T(\widehat{F}x_{k}-\widehat{F}x_{q}, t)\leq 1-u \quad \textit{and} \quad T(\widehat{F}x_{k}-\mathcal{L}, t/2)>1-\varepsilon,
\end{eqnarray*}
in particular $T(\widehat{F}x_{q}-\mathcal{L}, t/2)>1-\varepsilon$. Then
\begin{eqnarray*}
1-u \geq T(\widehat{F}x_{k}-\widehat{F}x_{q}, t)\geq T(\widehat{F}x_{k}-\mathcal{L}, t/2)\bigstar T(\widehat{F}x_{q}-\mathcal{L}, t/2)>(1-\varepsilon)\bigstar(1-\varepsilon)>1-u,
\end{eqnarray*}
which is not possible. On the other hand,
\begin{eqnarray*}
U(\widehat{F}x_{k}-\widehat{F}x_{q}, t)\geq u  \quad \textit{and} \quad  U(\widehat{F}x_{k}-\mathcal{L}, t/2)<\varepsilon,
\end{eqnarray*}
in particular $U(\widehat{F}x_{q}-\mathcal{L}, t/2)<\varepsilon$. Then
\begin{eqnarray*}
u\leq U(\widehat{F}x_{k}-\widehat{F}x_{q}, t)\leq U(\widehat{F}x_{k}-\mathcal{L}, t/2)\bigstar U(\widehat{F}x_{q}-\mathcal{L}, t/2)< \varepsilon\blacktriangle\varepsilon<u,
\end{eqnarray*}
which is not possible. Hence $B(\varepsilon, t)\subset A(\varepsilon, t)$. Therefore, by \ref{equ:3} $d(B(\varepsilon,t))=0$. Hence,
$x$ is $FSCa-IFN$.
\end{proof}

\begin{defin}(\cite{Saad})
An IFNS $(X,T,U,\bigstar,\blacktriangle)$ is said to be complete if every Cauchy sequence
is convergent in $(X,T,U,\bigstar,\blacktriangle)$.
\end{defin}

\begin{defin}
An IFNS $(X,T,U,\bigstar,\blacktriangle)$ is said to be statistically $(FSC-IFN)$ complete is every
statistically ($FSC-IFN$, respectively) Cauchy sequence with respect to intuitionistic fuzzy norm $(T,U)$ is statistically
($FSC-IFN$, respecitvely) convergent with respect to to intuitionistic fuzzy norm $(T,U)$.

\end{defin}

\begin{thm}\label{theo:4}
An IFNS $(X,T,U,\bigstar,\blacktriangle)$ is $(FSC-IFN)-$complete.
\end{thm}

\begin{proof}
Let $x=(x_{k})$ be $(FSC-IFN)-$Cauchy but not $(FSC-IFN)-$convergent. For a given $\varepsilon>0$
and $t>0$, choose $u>0$ such that $(1-\varepsilon)\bigstar(1-\varepsilon)>1-u$ and $\varepsilon\blacktriangle\varepsilon<u$.
Now
\begin{eqnarray*}
T(\widehat{F}x_{k}-\widehat{F}x_{N}, t)\geq T(\widehat{F}x_{k}-\mathcal{L}, t/2)\bigstar T(\widehat{F}x_{N}-\mathcal{L}, t/2)>(1-\varepsilon)\bigstar(1-\varepsilon)>1-u
\end{eqnarray*}
and
\begin{eqnarray*}
U(\widehat{F}x_{k}-\widehat{F}x_{N}, t)\leq U(\widehat{F}x_{k}-\mathcal{L}, t/2)\blacktriangle U(\widehat{F}x_{N}-\mathcal{L}, t/2)<\varepsilon \blacktriangle \varepsilon<u,
\end{eqnarray*}
since $x$ is not $(FSC-IFN)-$convergent. Therefore $d(E^{C}(\varepsilon, t))=0$, where
\begin{eqnarray*}
E(\varepsilon,t)=\{k\in \mathbb{N}: U_{x_{k}-x_{N}}(\varepsilon)\leq 1-r\}
\end{eqnarray*}
and so $d(E(\varepsilon,t))=1$, which is a contradiction, since $x$ was $(FSC-IFN)-$Cauchy. So that $x$ must be
$(FSC-IFN)-$convergent. Hence every IFNS is $(FSC-IFN)-$complete.
\end{proof}

We can give following theorem from Theorems \ref{theo:2}, \ref{theo:3}, \ref{theo:4}.

\begin{thm}
Let $(X,T,U,\bigstar,\blacktriangle)$ be an IFNS. Then, for any sequence $x=(x_{k})\in X$,
the following conditions are equivalent:\\
\begin{itemize}
\item [i.] $x$ is $(FSC-IFN)-$convergent.\\
\item[ii.] $x$ is $(FSC-IFN)-$Cauchy.\\
\item[iii.] IFNS $(X,T,U,\bigstar,\blacktriangle)$ is $(FSC-IFN)-$complete.\\
\item[iv.] There exists an increasing index sequence $J=(j_{n})$ of natural numbers such that
$d(K)=1$ and the subsequence $(x_{k_{n}})$ is a $(FSC-IFN)-$Cauchy.
 \end{itemize}
\end{thm}

\section{Conclusion}

In this paper, the well-known concept of statistical convergence of real numbers combine
 with the definition of Fibonacci numbers is studied in intuitionistic fuzzy normed spaces.
Firstly, the definition of Fibonacci statistical convergence with respect to intuitionistic fuzzy
normed space(FSC-IFN) is given and some properties of FSC-IFN are investigated. After,
the definitions of Fibonacci statistical Cauchy with respect to intuitionistic fuzzy
normed space(FSCa-IFN) and FSC-IFN completeness are given. It shows that an IFNS is (FSC-IFN)-complete.

\end{document}